\documentclass[8pt]{article}
\usepackage{amsmath,amsfonts,amssymb,amsthm}
\usepackage{hyperref}
\usepackage{graphicx}
\usepackage{tikz} 
\usepackage[utf8]{inputenc}

\newtheorem{thm}{Theorem}

\newtheorem{prop}{Proposition}

\newcommand{\R}{{\mathbb R}}
\newcommand{\Z}{{\mathbb Z}}
\title{Quasistationary Distribution for the Invasion Model on a Complete Bipartite Graph} 
\author{Clayton Allard \and Iddo Ben-Ari\thanks{Corresponding author, {\tt iddo.ben-ari@uconn.edu}. Research performed in the 2021 and 2022 UCONN Markov Chains REU, partially funded by REU grants  H98230-20-1-0253 and H98230-21-1-0016 from the Mathematical Sciences Program in the NSA to Iddo Ben-Ari.}\and Shrikant Chand \and Van Hovenga \and Edith Lee \and Julia Shapiro} 
\date{}
\begin{document}
\maketitle
\begin{abstract} 
The Invasion Model on the complete bibartitle graph was introduced and studied by physicists as a rudimentary model for opinion dynamics on complex networks. We identify the limit of the Quasistationary distribution for the model as one partition size tends to infinity. The limit is a highly dispersed measure.  A distinctive feature of the model is that of two time scales with non-trivial interaction. The work and the results complement and are in sharp contrast to the analogous results on the closely related Voter Model. 

\end{abstract} 
\section{Introduction} 
\subsection{The Invasion and Voter Models} 
\label{sec:invasion} 
The model we study here, the Invasion Model (or Invasion Process),  was initially introduced in \cite{aip} as a ``reverse'' Voter Model, as part of a larger study of opinion dynamics on complex networks. Though both are usually considered in the continuous-time setting, here we chose to work with the discrete time versions (the skeleton processes corresponding to continuous-time Markov chains). This is also because our main  object of study, the Quasistationary Distribution, is identical for both continuous-time and discrete-time versions, provided that in the latter all exponential clocks assigned to  vertices all have the same rate. We begin with a description of the Invasion Model.  Let $G=(V,E)$ be a finite connected graph. We assume that the vertex set $V$ has at least two elements and that $G$ has no loops, that is, all edges contain two distinct vertices.  The Invasion Model is a discrete-time Markov chain ${\bf \eta} = \{\eta_t:t\in \Z_+\}$ on the set $\Omega=\{0,1\}^V$, the set of functions $\eta: V \to \{0,1\}$, where each function should be interpreted as an opinion assignment, with $\eta(v)$ representing the opinion of the ``agent'' at vertex $v$: $0$ for a  ``no'' and $1$ for a  ``yes''. Each such function will be referred to as an {\it opinion configuration}. Starting from an initial distribution on opinion configurations, that it a distribution for $\eta_0$,  each unit of time $t \in \Z_+$ we do the following, independently of the past:
\begin{enumerate} 
\label{item:v}
\item Uniformly sample a vertex $v$. 
\label{item:u}
\item Conditioned on the vertex  $v$, uniformly sample a neighbor $u$.
\item\label{item:invasion} Set $\eta_{t+1}(w) = \eta_t (w)$ for $w\ne u$ and set $\eta_{t+1}(u)=\eta_{t}(v)$. That is, $v$ imposes their opinion on $u$ (or ``invades'' to $u$). 
\end{enumerate} 
The Voter Model follows similar dynamics with the only difference being   \ref{item:invasion} replaced by 
\begin{enumerate} 
\item[3'.] Set $\eta_{t+1}(w) = \eta_t (w)$ for $w\ne v$ and set $\eta_{t+1}(v)=\eta_{t}(u)$. That is $v$, accepts the opinion on $u$.  
\end{enumerate}
Note that the dynamics for the Invasion and the Voter models are identical if and only if $G$ is a constant-degree graph  Both processes are Markov chains and the absorbing states are the constant functions, namely states of the form $\eta (\cdot)\equiv C$ for $C\in \{0,1\}$, the {\it consensus} states. We write $\Delta$ for the set of two consensus states.  
We study the Invasion model for the case where $G$ is the complete bi-partite graph $K_{m,n}$. That is, a graph with vertex set $V$ which is the disjoint union of $\mathcal{S}$ a ``small'' partition,  a set with $m\ge 1$ elements, and $\mathcal{L}$, a ``large'' partition with $n\ge m$ elements, and edge set $E$ consisting of  all pairs with one element from $\mathcal{S}$ and one element from $\mathcal{L}$. These graphs are probably the simplest to study except for the complete graphs, yet as our analysis suggests and our results show they exhibit non-trivial structure.  
\subsection{Some Motivation} 
Voter and Invasion (or more general opinion dynamics) on complete bi-partite graphs can be interpreted as rudimentary models for a social interaction  with a ``small'' class of highly connected agents and a ``large'' class of agents, each with a small set of connections. The Voter and the Invasion Model represent two extremely different cultures, a culture of following in the Voter Model and a culture of domination in the Invasion Model. Clearly, both will eventually reach a consensus. Intuitively, the Invasion Model is expected to reach a consensus at a slower pace than the Voter Model. To see this observe that on $K_{m,n}$ with $n\gg m$ the ordered pair $(v,u)$ sampled in steps 1.  and 2. above describing the Invasion dynamics is typically of the form $v \in \mathcal{L}$ and $u\in \mathcal{S}$. Thus, in the Invasion Model the large group will rarely see a change while the small group will ``flicker'', changing opinions rapidly. On the other hand, in the Voter model typically the small group will reach a consensus eventually carrying the large group with it. From these heuristics, one expects the ``pre-consensus'' distribution for the two models to be quite different: for the Voter Model it will be quite close to consensus while for the Invasion Model it will be very far from consensus. Pre-consensus distributions are often described by the Quasistationary Distributions (QSDs) or the essentially equivalent notion of Quasi-Limiting Distributions (QLDs), see Section \ref{sec:QSD_QLD}. The main result of our work, Theorem \ref{thm:main}, provides a complete description for the limit of the QSD for the Invasion Model on $K_{m,n}$  as $n\to\infty$. This work is related  to \cite{voter} where the analogous results for the Voter Model (extended to  $k$-partite graphs  in \cite{voterkpartite}), and the  differences will be discussed in Section \ref{sec:limitQSD}, and as the proofs show, either case has some distinct and distinctive mathematical features. Analysis of the Invasion Model reveals two different orders of magnitude which need to be considered at the same time. Analysis of the Voter Model reveals a power laws for the number of opinions dissenting from the eventual consensus opinion on the small group. The present work, as well as \cite{voter},  were initially motivated by results in physics literature \cite[Section III]{redner} \cite{Redner_3} \cite{redner_4} that among other topics looked at the expected time for consensus on  $K_{m,n}$ as a prelude for more complex networks, yet we were not able to find literature on QSDs for these models, except for the case of a complete graph (\cite[Section 5]{dickman}).

\subsection{A General Model, Survival Rate and QSD} 
\label{sec:general} 
In this section we consider a (slightly more) general opinion dynamics on a graph, with  Invasion and Voter models being special cases of. This more general perspective is not new and the structure (with some differences) has been referred to broadly as Voter Model \cite[Section 5]{Aldous_Interacting}, a practice which we will not follow here because for our main object of interest, the QSD, it is more convenient to clearly distinguish between what we previously defined as the Invasion Model  and the Voter Model. 

 To initiate the discussion, observe that by construction of both Invasion and Voter Model in Section \ref{sec:invasion} above, the flow of opinions, namely whose opinion is assigned to whom at each step, does not depend on the  opinion distribution but only on a mechanism of sampling an ordered pair of neighbors identifying the second as accepting the opinion of the first.  

As before, $G=(V,E)$ is a finite connected graph with vertex set $V$ having  at least two elements and no loops. Let $E_d = \{(v,u):\{u,v \}\in E\}$, the set of ordered pairs of distinct vertices in $V$, and let $\rho$ be a probability measure on $E_d$ with full support. The generalized process will be also denoted by ${\bf \eta}= (\eta_t:t\in \Z_+)$ taking its values in the state space of opinion configurations $\Omega$.  Starting from some (possibly random) initial opinion assignment $\eta_0$, at time  $t\in\Z_+$ sample $({\cal V},{\cal U})$ according to $\rho$, independently of the past, and set $\eta_{t+1}({\cal U}) =\eta_t ({\cal V})$, and $\eta_{t+1}(w)=\eta_t(w)$ for all $w\ne {\cal U}$. That is the opinion of ${\cal U}$ at time $t+1$ is set to be the opinion of ${\cal V}$ at time $t$. For all $u \ne {\cal U} $, $\eta_{t+1}(u) = \eta_{t}(u)$: the opinion remains unchanged from time $t$ to $t+1$. In what follows we will write $P_\gamma$ for the distribution of ${\bf \eta}$ with initial distribution (distribution of $\eta_0$) being $\gamma$. 

The resulting Markov chain has $\Delta$, the set of constant functions as its absorbing set. Due to the connectivity and finiteness of  $G$,  $\Delta$ is accessible from every state. In particular, letting 
$$\tau =\inf \{t \in \Z_+: \eta_t \in \Delta\},$$ 
we have that $\tau < \infty$ a.s. and $\tau$ has geometric tails (provided the support of the initial distribution is not contained in $\Delta$). The geometric tail of $\tau$ (which in general may depend on the initial distribution) is central to our work.  Let $S(\rho)$ denote the restriction of the transition function for ${\bf \eta}$ to $(\Omega-\Delta)\times (\Omega-\Delta)$. The spectral radius of $S(\rho)$, $\lambda_{OD}(\rho)$ (OD for ``Opinion Dynamics'') is in $(0,1)$. Let $\|\cdot\|_{\infty}$ denote the $\ell^\infty$-norm and the induced operator norm. By the  spectral radius formula $\lambda_{OD}(\rho)$  \begin{equation}
 \label{eq:spectral_radius} \lambda_{OD}(\rho) = \lim_{t\to\infty} \|S(\rho)^t\|_\infty^{1/t}. 
 \end{equation} 
Of course,
$$ \|S(\rho)^t \|_\infty = \|P_{\cdot} (\tau>t) \|_\infty=\max_{\mu} P_\mu (\tau>t),$$ 
(with maximum attained at a delta measure $\mu$),  and standard arguments then give
\begin{equation} 
\label{eq:log_tail}\ln \lambda_{OD}(\rho) = \lim_{t\to\infty} \frac{1}{t}\ln  \max_\mu  P_\mu(\tau>t) = \max_\mu \lim_{t\to\infty} \ln P_\mu (\tau>t) .
\end{equation} 
Therefore, $\lambda_{OD}(\rho)$ represents the (geometric) survival rate for the opinion dynamics.  The next section will be devoted to finding  an alternate expression for $\lambda_{OD}(\rho)$.

\subsection{QSD and QLD} 
\label{sec:QSD_QLD}
This section serves as a  quick reminder of the general notions of a QSD and QLD in the context of opinion dynamics. We refer the reader to \cite{qsdbook} for more details. 

Recall that a probability distribution  $\nu$ on the nonabsorbing set $\Omega-\Delta$ is called a QSD if 
\begin{equation*} 
 P_\nu (\eta_1 \in \cdot ~| \tau>1) = \nu(\cdot). 
\end{equation*} 
This is equivalent to $P_\nu (\eta_t \in \cdot~ |\tau>t) = \nu(\cdot)$ for all $t\in\Z_+$, as well as the eigenvalue problem 
\begin{equation} 
\label{eq:QSD_eigen} \nu S(\rho) = \lambda \nu,
\end{equation} 
for some $\lambda>0$. When $S(\rho)$ is irreducible, it follows from the Perron-Frobenius theorem that $\lambda$ is its spectral radius and that the solution is unique.  Similarly, the probability measure  $\nu$ on $\Omega-\Delta$ is called a Quasi-Limiting Distribution if there exist some probability measure $\mu$ on $\Omega-\Delta$ such that 
\begin{equation} 
\label{eq:QLD} \lim_{t\to\infty} P_\mu (\eta_t \in \cdot ~| \tau>t) = \nu(\cdot), 
\end{equation} 
where the limit is in distribution. It is well-known a $\nu$ is a QSD if and only if it is also a QLD.

The invasion model corresponds to  $\rho=\rho_I$ where 
\begin{equation} 
\label{eq:rhoI} \rho_I (v,u) = \frac{1}{|V|} \frac{{\bf 1}_{\{u,v\}\in E}}{\mbox{deg}(v)}.
\end{equation} 

The Voter model corresponds to  $\rho=\rho_V$ where 
\begin{equation} 
\label{eq:rhoV} \rho_V (u,v) = \rho_I (v,u).
\end{equation} 
Note that the two models coincide if and only if $G$ is a constant degree graph.  The simplest family of all non-constant degree graphs to study both models is the class of complete bi-partite graphs. 
\subsection{Reverse Flow, Coalescing Chains, and Survival Rate}
\label{sec:coalescing}
In this section we describe a process that is dual to ${\bf \eta}$ from last section. The duality is in the sense that the process here traces the flow of opinions backwards in time. The ideas of reverse dynamics and coalescing random walks or coalescing Markov chains presented here are not new and have been primarily used in the context of the Voter model, see for example  \cite{Aldous_Interacting}\cite{crw}\cite{Durrett_particles}\cite{Liggett}\cite{RIO_CRW}.  How is this tracing done? Recall that the opinion of vertex $u$ at time $t\ge 1$ is the opinion of ${\cal V}$  at time $t-1$ if $({\cal V},u)$ was sampled, that is if ${\cal U}=u$. Otherwise, $u$ keeps its opinion from time $t-1$.  These two  alternative can be viewed as for a particle moving from $u$ to ${\cal V}$ or staying put. Repeating the opinion tracing back in time,  due to the independence of the repeated samples $({\cal V},{\cal U})$, we see that 
\begin{itemize} 
\item Tracing back in time from any given vertex $u$ is a Markov chain (and therefore its dynamics can be extended indefinitely backward in time); and
\item The family of corresponding Markov chains indexed by $u\in V$ is on its own right a Markov chain.
\end{itemize} 
The latter object, the family of coupled Markov chains,  is the flow of opinions backward in time, which we now describe formally.

Let $\rho_2$ denote the marginal of the second component of $\rho$. That is 
$$ \rho_2 (u) = \sum_v \rho(v,u).$$ 
Clearly, $$ \rho (v,u) = \rho(v|u) \rho_2 (u).$$ 

Let ${\bf X}=(X_t(u):u\in V,t \in \Z_+)$ be the process defined as follows. 
\begin{enumerate}
    \item \label{eq:init} ${\bf X}_0(u)=u$ for all $u\in V$.   
    \item At each time $t\in \Z_+$ , sample ${\cal U}$ according to $\rho_2$, independently of the past and sample ${\cal V}$ according to $\rho(\cdot ~|~ {\cal U})$ (of course this is identical to sampling $({\cal V},{\cal U})$ according to $\rho$, it stresses  the time reversal). Set 
$$  X_{t+1}(u) = \begin{cases} {\cal V} & X_t (u) = {\cal U} \\ X_t (u) & \mbox{otherwise}. \end{cases} $$ 
\end{enumerate} 
Note that once $X_t(u)=X_t(u')$, then $X_{s}(u)=X_{s}(u')$ for all $s\ge t$, and that at each unit of time only processes at a single vertex will move. As $\rho$ has full support, these imply that eventually $X_t(\cdot)$ will be constant, namely all processes will be at the same place. When looking at each chain  $X_{\cdot} (u)$ individually, we observe that its transition function $\kappa$ is given by 
$$ \kappa (u,v) = \begin{cases} \rho(v,u)=\rho_2(u) \rho(v|u)  & v \ne u \\ 
1-\rho_2(u) & v=u.
\end{cases} $$ 

Thus, ${\bf X}$ is a system of coalescing Markov chains, and the analysis leading to the formal construction yields the following result: 
\begin{prop}
\label{prop:reverse}
Consider the opinion dynamics on a finite connected graph $G=(V,E)$ corresponding to the probability measure $\rho$ with full support on  $E_d=\{(u,v):\{u,v\}\in E\}$.  Let ${\bf X}$ be the coupled system defined above. Then for every $T\in \Z_+$, the process ${\bf X}$ restricted to the time interval $[0,T]$ has the same distribution as the backward flow of opinions. That is, the joint distribution of the random vectors $s\in [0,T], u\in V$ is the same as the joint distribution of the vertices whose opinion at time $T-s$ vertex $u$ holds, respectively. 
\end{prop} 

For distinct $u,v \in V$, let  $$\sigma_{u,v}=\inf\{t \in \Z_+: X_t(u)=X_t(v)\}.$$ 
Then by assumption on $\rho$, $\sigma_{x,y}<\infty$. Since it is an a.s. finite  stopping time for a finite-state chain, its distribution has  geometric tails. Let $\sigma = \sup_{u,v} \sigma_{u,v}$. Then, again, $\sigma$ is a.s. finite and its distribution has geometric tails. Moreover, the following holds
\begin{equation}\label{eq:lambdaCMC}  \lim_{t\to\infty} \frac{1}{t}  \ln P(\sigma >t)= \lim_{t\to\infty} \frac{1}{t} \ln \max_{u,v} P(\sigma_{u,v} >t) = \max_{u,v} \lim_{t\to\infty}\frac{1}{t} \ln P(\sigma_{u,v}>t),\end{equation} 
where all limits exist and is in $(-\infty,0)$. Write $\lambda_{CMC}(\rho)$ (CMC for ``Coalescing Markov Chains'') for the exponential of this limit. That is, $\lambda_{CMC}(\rho) \in (0,1)$ and it satisfies   
$$ P(\sigma>t) = \lambda_{CMC}(\rho)^{t (1+o(1))}.$$ 

We have the following 
\begin{prop}
\label{prop:samelambda}
 \begin{equation} \lambda_{CMC}(\rho) = \lambda_{OD}(\rho).
 \end{equation} 
\end{prop}
 This result reduces the computation of the tail of $\tau$ to the calculation of the ``heaviest'' tail among  all  coalescence times $\sigma_{v,v'}$ for the pairs of  weakly-coupled processes $({\bf X}_{\cdot}(v),{\bf X}_{\cdot}(v'))$.  

The proof of the Corollary we give here is essentially identical to the proof of \cite[Corollary 4.2]{voter}, with some change of notation\footnote{Specifically: $\lambda_{CRW}(G)$ there corresponds to  $\lambda_{CMC}(\rho)$ here and $\lambda_{V}(G,\mu)$ there is geometric tail of the distribution of $\tau$ under $P_\mu$, with $\max_{\mu} \lambda_V(G\mu)$ corresponding to our $\lambda_{OD}(\rho)$. Though the discussion there is on the voter model, the results generalize to our setting with the same proofs.}.    

\begin{proof}[Proof of Proposition \ref{prop:samelambda}] Observe that under any initial distribution, $\tau$ is stochastically dominated by $\sigma$ (whose distribution depends only on $\rho$). This translates into smaller tails and hence $\lambda_{OD}(\rho)\le \lambda_{CMC}(\rho)$. On the other hand, observe that the identification in Proposition \ref{prop:reverse} gives the following inequality. Fix distinct $(v,v')\in E_d$. If at time $t$, $(X_t(v),X_t(v'))=(u,u')$ for some $(u,u')\in E_d$ satisfying  $\eta_0(u)\ne \eta_0(u')$ then necessarily $\tau>t$.  Putting this together we have 
\begin{align*} 
P_\mu(\tau>t) &\ge \sum_{(u,u')\in E_d} P (({X}_t(v),{ X}_t(v'))=(u,u'),\eta_0(u)\ne \eta_0(u'))\\
& \ge P(\sigma_{v,v'}>t) \min_{(u,u')\in E_d} \mu(\eta_0(u)\le\eta_0(u')).
\end{align*} 
In particular, choosing $\bar \mu$ to be uniform over all $V$ configurations with exactly one vertex having opinion $1$ and all  the rest having opinion $0$, $\bar \mu(\eta_0(u) \ne \eta_0(u'))= \frac{2}{|V|}$ for every $(u,u')\in E_d$ so we conclude that 
$$ P_{\bar \mu} (\tau>t) \ge \frac{2}{|V|} P(\sigma_{v,v'}>t),$$ 
and as the lefthand side is valid for every $(v,v')\in E_d$, it follows from \eqref{eq:lambdaCMC} that $\lambda_{OD}(\rho) \ge \lambda_{CMC}(\rho)$. Since the reverse inequality also holds, the proof is complete.
\end{proof} 

We comment that both Proposition \ref{prop:reverse} and Proposition \ref{prop:samelambda} are routine extensions to well-known (but not extensively documented) analogous results for the Voter model.

Returning to the Invasion and Voter model, from \eqref{eq:rhoI} and \eqref{eq:rhoV}, we obtain the following respective second marginals $\rho_{I,2}$, $\rho_{V,2}$ and transition functions $\kappa_I$ and $\kappa_V$.

$$ \rho_{I,2}(u)= \frac{1}{|V|}\sum_{\{v:\{u,v\}\in E\}} \frac{1}{\mbox{deg}(v)},\quad \rho_{V,2}(u) = \frac{1}{|V|}.$$ 

For $u\ne v$:  
$$ \kappa_I  (u,v) = \frac{1}{|V|}\frac{{\bf 1}_{\{u,v\}\in E}}{\mbox{deg}(v)},\quad \kappa_V (u,v) = \frac{1}{|V|}\frac{{\bf 1}_{\{u,v\}\in E}}{\mbox{deg}(u)}.$$ 

Note that $\kappa_V$ corresponds to a  lazy random walk: with probability $1-\frac{1}{|V|}$, the process stays put, and otherwise it moves uniformly to one of the neighbors. On the other hand (with the exception of a constant-degree graphs), $\kappa_I$ corresponds to a Markov chain which moves to neighbors with probabilities proportional to the reciprocal to the degrees. 



\section{QSD for the Invasion on $K_{m,n}$}
In this section we collect the observations and results specific to the Invasion Model on $K_{m,n}$. 
We begin by introducing an induced chain on a ``simpler'' state space which will be used in our discussion of the model.

\subsection{Induced Chain and QSD } 
\label{sec:induced} 
As our main object of interest is the invasion process on $K_{m,n}$, in the remaining of this section we will consider this particular graph, but still under the general dynamics we introduced and discussed in Sections \ref{sec:general} and \ref{sec:coalescing}. 

Let $\Gamma$ denote the set of two states in $\Omega$ corresponding to a consensus opinion on each of ${\cal S}$ and ${\cal L}$,  with the opinions on each partitions being different. Then $\Gamma$ is inaccessible from any state. This implies that any left  eigenvector for $S(\rho)$ corresponding to a nonzero eigenvalue has zero for all entries in $\Delta$. This implies that a QSD is supported on the $\Omega-(\Delta\cup \Gamma)$, and that the spectral radius of $S(\rho)$ (on $\Omega- \Delta$) and of its restriction to $S(\rho)-(\Omega-\Delta)$ coincide.  To eliminate trivialities,  we will always assume 
 \begin{equation}
\label{eq:irreducible_condition} 
2\le m \le n\mbox{, or }~  m=1 \mbox{ and }n\ge 3.
\end{equation} 
A straightforward argument  then guarantees that the restriction of $S(\rho)$ to $\Omega- (\Delta \cup \Gamma)$ is irreducible. 

In what follows, we will somewhat abuse notation and consider $S(\rho)$ as the substochastic matrix on $\Omega- (\Delta\cup \Gamma)$. As noted above, this has no impact on the QSD or on the spectral radius. Our assumption \eqref{eq:irreducible_condition} yields the irreducibility and allows to apply \cite[Theorem 2.2]{voter} to conclude that the Invasion model on $K_{m,n}$ has a unique QSD which is also a QLD, and which is a left  eigenvector for $S(\rho_I)$ corresponding to eigenvalue $\lambda_{OD}(\rho_I)$, the spectral radius of $S(\rho)$. According to proposition \ref{prop:samelambda}, $\lambda_{OD}(\rho_I)=\lambda_{CMC}(\rho_I)$. We will denote the QSD by $\nu_{m,n}$. As a first step towards analyzing $\nu_{m,n}$ we will estimate (calculate in the special case $m=1$) the latter in Proposition \ref{prop:lambdacmc}. 

 Because the dynamics is invariant under permutation of the vertices within each of the partitions (namely relabeling of vertices in ${\cal S}$ or of vertices in ${\cal L}$), the uniqueness of the  QSD  implies that the QSD assigns the same probability to any two states  obtained from each other  through such permutations. As every two states $\eta,\eta'$ with the property that 
 $$\begin{cases}  \sum_{u \in {\cal S}} \eta (u) = \sum_{u\in {\cal S}} \eta' (u) & \mbox{ and } \\
 \sum_{u \in {\cal L}} \eta (u) = \sum_{u\in {\cal L}} \eta' (u).
 \end{cases}$$ 
  can be obtained from each other through such permutations, it readily follows that the QSD is a function of the numbers of ``yes'' in ${\cal S}$ and in ${\cal L}$. The invariance of the dynamics under these permutations yields an induced Markov chain on the state space of ordered pairs  $\{0,\dots,m\}\times\{0,\dots,n\}$, where an element  $(k,l)$ can be identified with the set 
  \begin{equation} 
  \label{eq:pair_eta} 
  \{\eta \in \Omega : \sum_{u \in {\cal S}} \eta (u)=k\mbox{ and }\sum_{u \in {\cal L}} \eta(u)=l\}.
  \end{equation} 
That is $(k,l)$   corresponds representing all $\binom{m}{k} \times \binom{n}{l}$ states in $\Omega$ with $k$ ``yes'' in  ${\cal S}$ and $l$  ``yes'' in ${\cal L}$. The induced Markov chain inherits its properties from the population dynamics. In particular, 
\begin{itemize} 
\item The absorbing set for the induced chain is the image of $\Delta$, $\{(0,0),(m,n)\}$; 
\item The set of inaccessible states for the induced chain is the image of $\Gamma$, $\{(0,n),(m,0)\}$;
\item All remaining states form an irreducible class;
\item The induced chain has a unique QSD supported on the remaining states the last item which is also a QLD, and for simplicity we will also denote by $\nu_{m,n}$, because each state  $(k,l)$ of the induced chain is a set of states in $\Omega$.  Clearly,  
\begin{equation}
\nu_{m,n}(k,l) = \binom{m}{k}\binom{n}{l}  \nu_{m,n} (\eta), 
\end{equation} 
where $\eta$ is any state in the set \eqref{eq:pair_eta}, that is a state with $k$ ``yes'' in ${\cal S}$ and $l$ ``yes'' in ${\cal L}$. 
\end{itemize}

\subsection{Transitions for Induced Chain}
\label{sec:transitions}
To derive the transition probabilities for the induced chain, we first observe that for the Invasion Model on $K_{m,n}$, the kernel $\rho_I$ \eqref{eq:rhoI} is given by 
\begin{equation} 
\label{eq:rho_Kmn}\rho_I(v,u) = \frac{1}{m+n} \begin{cases} \frac{1}{m} & (v,u) \in {\cal L} \times {\cal S}
\\ \frac{1}{n} & (v,u) \in {\cal S} \times {\cal L}.\end{cases}
\end{equation} 
Therefore the  transitions for the corresponding induced chain are: 
\begin{enumerate}
    \item $(k, l) \rightarrow (k+1,l): \frac{l}{n+m}\cdot \frac{m-k}{m} = \frac{l(m-k)}{m(n+m)}.$ First choose one of the $l$ vertices in group $\mathcal{L}$ with opinion $1$ then choose one of the $m-k$ vertices in group $\mathcal{S}$ with opinion $0$ to convert into opinion $1.$
    \item $(k, l) \rightarrow (k-1,l): \frac{n - l}{n+m} \cdot \frac{k}{m} = \frac{k(n-l)}{m(n+m)}.$ First choose one of the $n - l$ vertices in group $\mathcal{L}$ with opinion $0$ then choose one of the $k$ vertices in group $\mathcal{S}$ with opinion $1$ to convert into opinion $0.$
    \item $(k, l) \rightarrow (k,l+1): \frac{k}{n+m} \cdot \frac{n-l}{n} = \frac{k(n-l)}{n(n+m)}.$ First choose one of the $k$ vertices in group $\mathcal{S}$ with opinion $1$ then choose one of the $n-l$ vertices in group $\mathcal{L}$ with opinion $0$ to convert into opinion $1.$
    \item $(k, l) \rightarrow (k,l-1): \frac{m-k}{n+m} \cdot \frac{l}{n} = \frac{l(m-k)}{n(n+m)}.$ First choose one of the $m-k$ vertices in group $\mathcal{S}$ with opinion $0$ then choose one of the $l$ vertices in group $\mathcal{L}$ with opinion $1$ to convert into opinion $0.$
    \item $(k, l) \rightarrow (k,l): \frac{k}{n+m}\cdot \frac{l}{n} + \frac{m-k}{n+m}\cdot\frac{n-l}{n} + \frac{l}{n+m}\cdot \frac{k}{m} + \frac{n-l}{n+m}\cdot\frac{m-k}{m} = \frac{kl + (m-k)(n-l)}{nm}.$ First choose any vertex and then choose a neighbor with the same opinion.
\end{enumerate}

\subsection{Reverse Flow and Survival Rate}
As noted in Section \ref{sec:induced}, the unique QSD for the Invasion Model on $K_{m,n}$ is a left eigenvector for $S(\rho)$ corresponding to its spectral radius $\lambda_{OD}(\rho_I)$. By Proposition \ref{prop:samelambda}, $\lambda_{OD}(\rho_I)=\lambda_{CMC}(\rho_I)$, and the goal of this section is to exploit the relatively simple structure of the coalescing Markov chains to identify this quantity. The first step is to describe the reverse flow. The second step is to look at the structure of pairs of chains in the reverse flow and their coalescence. The last and third step is to derive a formula for the  geometric tail of the coalescence time, $\lambda_{CMC}(\rho_I)$. 

\subsubsection*{The Reverse Flow} 
In order to determine the reverse flow, we first record the marginal and conditional distributions for the Invasion Model on $K_{m,n}$ from the formula for $\rho_I$,   \eqref{eq:rho_Kmn}. Let $s\in {\cal S}$ and let $\ell \in {\cal L}$. Then 
\begin{align*}
\rho_{I,2}(s)& = \frac{n}{m+n}\frac{1}{m} &\rho_{I,2}(\ell) = \frac{m}{m+n}\frac{1}{n}  \\
\rho_I(\ell|s) &= \frac{1}{n}& \rho_I(s|\ell) = \frac{1}{m}.
\end{align*} 
The resulting reverse flow is then: 
\begin{itemize} 
\item Sample a vertex ${\cal U}$ according to $\rho_{I,2}$, which in this case assigns probability $\frac{n}{m+n}\frac{1}{m}$ to each vertex in ${\cal S}$ and probability $\frac{m}{m+n}\frac{1}{n}$ to each vertex in ${\cal L}$. 
\item Sample a neighbor ${\cal V}$ uniformly. 
\item Move all chains at ${\cal U}$ to ${\cal V}$ and have all other chains stay put. 
\end{itemize} 

\subsubsection*{Pairs of Reverse Chains} 
\begin{figure}[t]
    \centering
    \resizebox{0.6\textwidth}{0.6\textwidth}{%
     \begin{tikzpicture}
\node (a1) at (0,0)  {};
\node (a2) at (4,0)   {};
\node (a3) at (2,-4)  {};

 \node at (-0.2,0.2) {$1$};
 \node at (3.8,0.2) {$2$};
 \node at (1.8,-3.8) {$3$};
 node 

\draw[rounded corners] (a1) rectangle (2.5,-2.5)  {};
\draw[rounded corners] (a2) rectangle (6.5,-2.5) {};
\draw[rounded corners] (a3) rectangle (4.5,-6.5) {}; 
\draw (a1)+(0.75,-1.25) ellipse (0.4 and 0.8) {}; 
\draw (a1)+(1.75,-1.25) ellipse (0.3 and 0.6) {};

\draw (a2)+(0.75,-1.25) ellipse (0.4 and 0.8) {}; 
\draw (a2)+(1.75,-1.25) ellipse (0.3 and 0.6) {};

\draw (a3)+(0.75,-1.25) ellipse (0.4 and 0.8) {}; 
\draw (a3)+(1.75,-1.25) ellipse (0.3 and 0.6) {};

\draw[fill=black]  (a1)+(0.75,-1.25-0.4) circle (0.08) {}; 
\draw[fill=black]  (a1)+(0.75,-1.25+0.4) circle (0.08) {};

\draw[fill=black]  (a2)+(1.75,-1.25-0.4) circle (0.08) {}; 
\draw[fill=black]  (a2)+(1.75,-1.25+0.4) circle (0.08) {};

\draw[fill=black]  (a3)+(0.75,-1.25) circle (0.08) {}; 
\draw[fill=black]  (a3)+(1.75,-1.25) circle (0.08) {}; 

\draw[<->] (a1)+(1.25,-2.6)  -- (3.2,-3.9) node[midway,anchor=east] {};
\draw[<->] (a2)+(1.25,-2.6)  -- (3.6,-3.9) node[midway,anchor=west] {};

\draw[->] (a3)+(1.25,-2.6) -- (3.25,-8.6) node[anchor=north]  {Coalesce};

\node (a1s) at (0,-1.25) {} ;
\draw [->] (-0.1,-1) to [out=120,in=210,looseness=5] (-0.1,-1.5) {};
\draw [->] (6.5+0.1,-1) to [out=60,in=-30,looseness=5] (6.5+0.1,-1.5) {};
\draw [->] (2-0.1,-5) to [out=120,in=210,looseness=5] (2-0.1,-5.5) {};
\end{tikzpicture} 
}
    \caption{The transition diagram for the coalescing Markov chains. Top left box: both chains  are in ${\cal L}$ (represented by large ellipse), case 1 above. Top right box: both chains are  in ${\cal S}$ (represented by small ellipse), case 2  above. Center block: one chain in ${\cal L}$ and one in ${\cal S}$, case 3 above.}
 \label{fig:chains}
\end{figure}
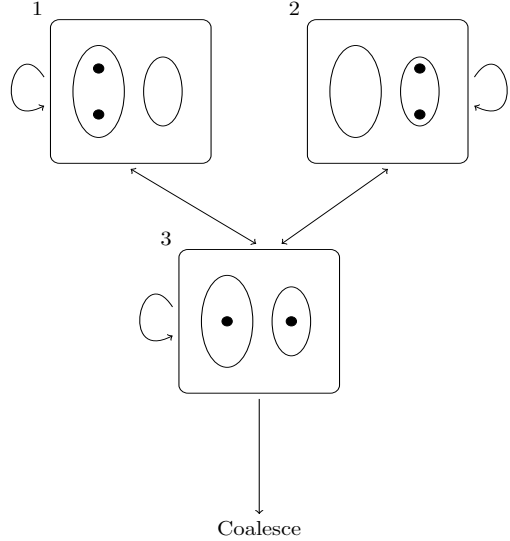

Fix distinct vertices $u,u'\in K_{m,n}$ and examine the behavior of $({  X}_t(u),{   X}_t(u'):t\in\Z_+)$ restricted to the time interval before the two processes coalesce.  Initially,  ${   X}_0(u)=u$ and ${   X}_0(u')=u'$. Before they coalesce, each unit of time $t$, the pair of Markov chains move across the three states below. The transition diagram is illustrated in Figure \ref{fig:chains}. 

\begin{enumerate} 
\item
\label{case:bothL} ${   X}(u),{   X}(u')\in {\cal L}$ (with ${   X}(u)\ne{   X}(u'))$. Then 
\begin{enumerate}
    \item With probability $2\rho_2(\ell) =\frac{2}{n} \frac{m}{m+n}$ one of the two will move to $\mathcal{S}$ and the other will stay put. This gets us to case \ref{case:diffs}.
    \item With the remaining probability both will stay put.
 \end{enumerate}
 \item
 \label{case:bothS} ${   X}(u),{   X}(u')\in {\cal S}$ (with ${   X}(u)\ne{   X}(u'))$. Then: 
\begin{enumerate} 
\item With probability $2\rho_2(s) = \frac{2n}{m+n}\frac{1}{n}$ one of the two will move to $\mathcal{L}$ and the other will stay put. This gets us to state \ref{case:diffs}. 
\item With the remaining probability both will stay put. 
\end{enumerate} 
\item \label{case:diffs}  ${   X}(u)\in {\cal S}$ and ${   X}(u') \in{\cal  L}$, or ${   X}(u) \in {\cal L} $ and ${   X}(u')\in {\cal S}$. 
\begin{enumerate} 
\item With probability $\rho_2(s) =\frac{n}{(m+n)m}$ the chain in $\mathcal{S}$ moves to $\mathcal{L}$, where with  probability $\frac{1}{n}$, the two chains coalesce and with the remaining probability, $\frac{n-1}{n}$,  we are back in case \ref{case:bothL}.
\item With probability $\rho_2(\ell) =\frac{m}{(m+n)n}$, the chain in $\mathcal{L}$ moves to $\mathcal{L}$ where with probability $\frac{1}{m}$ the two chains coalesce and with the remaining probability, $\frac{m-1}{m}$, we are back in case \ref{case:bothS}. 
\item With probability $1-\rho_2(s)-\rho_2(\ell) =1 -\frac{m^2+n^2}{(m+n)mn}$ both chains stay put. 
\end{enumerate}  
\end{enumerate}

\subsubsection*{Survival Rate}
These dynamics are Markovian, as they depend only on the current state of the system (one of the three listed above) and the corresponding substochastic matrix $p$ is
\begin{equation}
p=
\begin{pmatrix}\label{eq:tm}
1 - \frac{2m}{m+n}\frac{1}{n}  & 0 & \frac{2m}{m+n}\frac{1}{n} \\
0 & 1- \frac{2n}{m+n}\frac1{m} & \frac{2n}{m+n}\frac{1}{m}\\
\frac{n-1}{m+n}\frac{1}{m} & \frac{m-1}{m+n}\frac{1}{n} & 1- \frac{m^2+n^2}{(m+n)mn}.
\end{pmatrix}
\end{equation}
Note that $\lambda_{CMC}(\rho_I)$ is the spectral radius for $p$ by definition. We have the following: 
\begin{prop} 
\label{prop:lambdacmc} 
Consider the Invasion model on $K_{m,n}$ with $m,n$ satisfying \eqref{eq:irreducible_condition}. Then 
\begin{enumerate} 
\item 
\begin{equation}\label{eq:lambdaI}\lambda_{CMC}(\rho_I) = 1 - \frac{2m}{(m+n)n^2}+ o(n^{-3}), 
\end{equation} 
as $n\to\infty$. 
\item If $m=1$ 
\begin{equation}
\begin{split}
\label{eq:lambdaI1}
    \lambda_{CMC}(\rho_I)   &=1- \frac{ (3+n^2) - \sqrt{(3+n^2)^2 - 8(n+1)}}{2n(n+1)}\\
     & =1-\frac{2}{(3+n^2)n} + \Theta(n^{-6}),
\end{split} 
\end{equation} 
as $n\to\infty$.
\end{enumerate} 
\end{prop}
The analogous quantity for the Voter Model on $K_{m,n}$ is
\cite[Proposition 5.1]{voter}. 
$$ \lambda_{CMC}(\rho_V) = 1 - \frac{2}{m+n}\left ( 1 -\sqrt{1- \frac{1}{2m}-\frac{1}{2n}}\right).$$ 
The analysis for the Voter Model on $K_{m,n}$ is simpler because the reverse chains are random walks, allowing to reduce the matrix analogous to $p$ from \eqref{eq:tm} to a $2\times 2$ matrix. As under the QSD the absorption time is geometric with parameter $1-\lambda_{CMC}$, we observe that the corresponding expectation is $\Theta(n^3)$ for the  Invasion Model and $\Theta(n)$ for the Voter Model. The source for the difference is that the Invasion dynamics make it harder to reach consensus on ${\cal L}$. The gap between the behavior of the two models is further magnified when considering the respective QSDs, the topic of the next section.  
\begin{proof}[Proof of Proposition \ref{prop:lambdacmc}]
Let $\lambda$ denote the spectral radius (or Perron eigenvalue) for $p$. We  subtract $1$ from the diagonal  of $ p$ and multiply the matrix by $(m+n)nm$. This yields the matrix 
\begin{equation} 
\label{eq:pbarmatrix} \bar p = \begin{pmatrix} -2 m^2 & 0 & 2 m^2 \\ 
0 & -2n^2 & 2n^2 \\ 
(n-1) n & (m-1)m & - m^2 - n^2 
\end{pmatrix} 
\end{equation} 

The Perron eigenvalue for $\bar p$ is therefore  $\bar \lambda = (\lambda -1)(m+n)nm$. Note that if $m=1$, the second state is inaccessible, reducing the problem to a simpler two-state chain problem. We will first treat the more general case $m\ge 2$ and return to the special case $m=1$ at the end. 

We will consider $\bar p$ as the generator of a continuous-time Markov chain on the states $1,2,3$ and an absorbing state $4$, accessible only from state $3$. Let $f_1,f_2,f_3$ denote the expected absorption time from each of the respective states. Then we have 

$$ f_1 = \frac{1}{2m^2} + f_3,f_2 = \frac{1}{2n^2} + f_3$$ 
and 
$$ f_3 = \frac{1}{m^2+n^2} + \frac{(n-1)n}{m^2+ n^2}f_1 + \frac{m(m-1)}{m^2+n^2}f_2 + \frac{m+n}{m^2+n^2}.$$ 

This gives 

$$ f_3 (m+n)  = \left ( 1 + \frac{(n-1)n}{2m^2} + \frac{(m-1)m}{2n^2} + m+n \right).$$ 
Therefore, $f_3 \sim \frac{n}{2m^2}$, and as a result, $f_1,f_2\sim f_3$. If $\bar \pi$ represents the QSD for $\bar p$, then $\bar \pi \bar p= \bar \lambda \bar \pi$, and in particular under  $\bar \pi$, the absorption time is exponential with parameter $-\bar \lambda$. Thus,  the expected absorption time under $\bar \pi$ is $-\frac{1}{\bar \lambda}$. Since $f_1\sim f_2\sim f_3\sim \frac{n}{2m^2}$, it follows that the expected absorption time from any initial distribution is $\sim \frac{n}{2m^2}$. This gives $\bar \lambda \sim - \frac{2m^2}{n}$. From the definition of $\bar \lambda$, we finally obtain $\lambda -1\sim -\frac{2m}{(m+n)n^2}$. As $ \frac{1}{(m+n)^2 n}-\frac{1}{(m+n)n^2}=\Theta(n^{-4})$, giving  \eqref{eq:lambdaI}.

It remains to prove \eqref{eq:lambdaI1}. Remove state 2 from the matrix in \eqref{eq:pbarmatrix} we are left with

\begin{equation*} 
 \bar{\bar p} = \begin{pmatrix} -2  &  2  \\ 
(n-1) n &  - 1 - n^2 
\end{pmatrix} 
\end{equation*}
The characteristic polynomial is  $\lambda^2 +(3+n^2) \lambda +2(n+1),$ and therefore the Perron eigenvalue is  
$$  \bar {\bar \lambda}= \frac{ -(3+n^2) + \sqrt{(3+n^2)^2 - 8(n+1)}}{2},$$ 
and so 
$$ \lambda =1- \frac{ (3+n^2) - \sqrt{(3+n^2)^2 - 8(n+1)}}{2n(n+1)}. $$ 

\end{proof} 
\subsection{Limit for the QSD} 
\label{sec:limitQSD}
In this section we present and prove our main result on the limit of the QSD for the Invasion Model on $K_{m,n}$ as $n\to \infty$.  Recall that $\nu_{m,n}$ is the unique QSD for the Invasion Model on $K_{m,n}$, where $m$ and $n$ satisfy \eqref{eq:irreducible_condition}. We showed in Section \ref{sec:induced}, that $\nu_{m,n}$ is a function of the number of ``yes'' in ${\cal S}$ and the number of ``yes'' in ${\cal L}$ and observed that $\nu_{m,n}$ can be therefore viewed as a function on $\{0,\dots,m\}\times \{0,\dots,n\}$ where $\nu_{m,n}(k,l)$ is the probability $\nu_{m,n}$ assigns to the set of $\binom{m}{k}\times \binom{n}{l}$ states with $k$ ``yes'' in ${\cal S}$ and $l$ ``yes'' in ${\cal L}$. As we will let $n\to\infty$, we will replace the number of ``yes'' opinions in ${\cal L}$ with the corresponding density of ``yes'' opinions in ${\cal L}$, allowing to treat the second marginal as a probability measure on $[0,1]$. We write $\bar \nu_{m,n}$ for the resulting probability measure. That is 
\begin{equation}
\label{eq:QSD} 
\bar \nu_{m,n} (k,dx) =  \nu_{m,n} (k,x n ) \delta(dx)_{\{0,\frac{1}{n},\dots,1\}}.
\end{equation} 
Equivalently, for any continuous $f: \{0,\dots,m\}\times [0,1]\to {\mathbb R}$, 
$$ \int f (k,x) d \bar \nu_{m,n}(k,x) = \sum_{k=0}^m \sum_{l=0}^n f(k,l/n) \nu_{m,n}(k,l).$$ 

 Here is our main result: 
\begin{thm}
\label{thm:main} 
$$\bar \nu_{m,n}(k,dx)  \underset{n\to\infty}{\Rightarrow} \binom{m}{k} x^k (1-x)^{m-k} dx,$$ 
where ``$\Rightarrow$'' denotes convergence in distribution and ``$dx$'' denotes the Lebesgue measure.   In particular, under the limit  distribution 
\begin{enumerate}
\item Each of the marginals is uniform. 
\item \label{item:first_marginal}The first marginal conditioned on the second $=x$ is ${\bf Bin}(m,x)$.
\item 
\label{item:second_marginal} The second marginal conditioned on the first $=k$ is $\mbox{\bf Beta}(k+1,m+k+1)$. 
\end{enumerate} 
\end{thm} 
We wish to make a number of comments. 

The statement in item \ref{item:first_marginal} is intuitive and  easier to prove. This is because the small group is constantly invaded by the large group, with  invading and invaded vertices sampled uniformly. Thus, the opinions of members in the small group should exhibit independent opinions with probability of ``yes'' equal to the density of ``yes'' in the large group, whatever the latter is. The changes in the large group are more rare and so harder to be traced. The distribution of the large group corresponds to terms of smaller order in the eigenvalue equation for the QSD. Though the entire eigenvalue equation tends to a triviality (this is because the corresponding eigenvalue $\lambda_{OD}(\rho_I)$ tends to $1$), the terms associated primarily with the small group in the eigenvalue equation for the QSD are of larger order than those associated to the large group. Decoupling the two and analyzing the latter requires refined analysis, even merely from the perspective of obtaining sufficient number of terms in the expansion of the eigenvalue $\lambda_{OD}(\rho_I)$. Indeed, most of the effort in our proof is in identifying the asymptotic distribution of the large group.

The QSD for the Voter Model on $K_{m,n}$ converges as $n\to\infty$ to a very different limit \cite[Theorem 1.1]{voter}: all are in consensus, except for finitely many dissenting vertices in ${\cal L}$. The consensus opinion is ``yes'' with probability $\frac 12$ and the number of dissenting vertices follows a Sibuya distribution, a distribution on the natural numbers with a power-law tail depending on $m$. Thus, not only the Invasion Model typically takes more time to reach consensus (discussion below Proposition \ref{prop:lambdacmc}), but this is also done along quite a different path. While the  pre-consensus structure of the Voter Model is very close to consensus,  the pre-consensus structure of the Invasion Model is quite the opposite being very ``disorganized'' with opinions of the small group essentially being IID. This seems to suggests that in the Invasion Model consensus is reached ``abruptly'' rather than ``gradually'', reminiscent of the cut-off phenomenon and a potential research direction we have not explored. 
\begin{figure}[ht!]
    \label{fig:QSD 3D graph}
    \centering
    \includegraphics[width = 10cm]{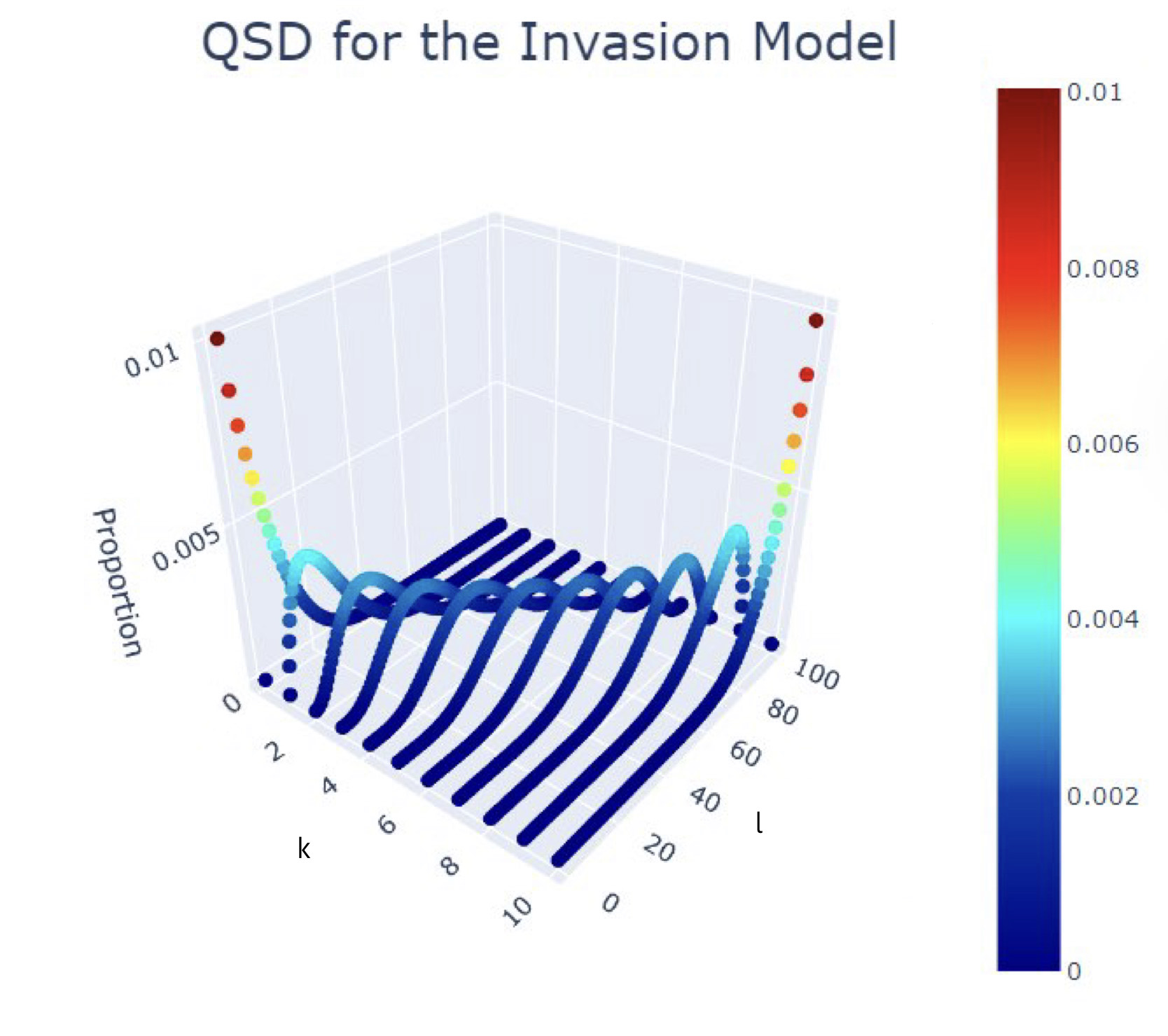}
    \caption{Simulated QSD for the Invasion Model with  on $K_{10,100}$. Each ``curve'' is a rendering of the mass function along a fixed value of $k$ (number of ``yes'' in ${\cal S}$) for $l$ (number of ``yes'' in ${\cal L}$) ranging from $0$ (bottom) to $100$, which (when scaled and normalized) is an approximation to a corresponding Beta distribution in  Theorem \ref{thm:main}, part \ref{item:second_marginal}.}
\end{figure}
\subsection{Proof of Theorem \ref{thm:main}}
We will use the fact that $\bar \nu_{m,n}$ are probability measures on a compact space $\{0,\dots,m\}\times [0,1]$. This implies that every subsequence has a convergent sub-subsequence. Since convergence in distribution is metric, it is enough to show that the limits along these sub-subsequences are all the same. 

To ease notation, we omit the dependence on $m$ and $n$ and write  $\nu, \bar \nu$ for $\nu_{m,n},\bar \nu_{m,n}$ respectively. 
\subsubsection{Rearrangement} 
First write the eigenvalue equation for the QSD. Recall that the absorbing set $\Delta$ is $\{(0,0),(m,n)\}$, and that $\nu$ satisfies the  equation $\lambda_{OD}(\rho_I) \nu = \nu S(\rho_I) $. Using the transition scheme in Section \ref{sec:transitions}, this equation becomes 
 \begin{align}
\label{eq:small1}
\nonumber
\lambda \nu (k,l)  &=\\
&\quad\frac{l}{(m+n)m}\left (  \nu(k-1,l) (m-k+1)+\nu(k,l)  k\right) \\
\label{eq:small2}
&+  \frac{n-l}{(m+n)m} \left (\nu (k+1,l) (k+1)+\nu(k,l) (m-k) \right) \\
\label{eq:large1}
& + \frac{k}{(m+n)n} \left( \nu(k,l-1) (n-l+1) + \nu(k,l) l\right)\\
\label{eq:large2}
& + \frac{m-k}{(m+n)n}\left( \nu(k,l+1) (l+1) + \nu(k,l) (n-l)\right),
\end{align} 
for any $(k,l) \in\{0,\dots,m\}\times\{0,\dots,n\}-\Delta$ (there is no need to eliminate the unattainable states $(0,n)$ and $(m,0)$ because in each of these cases the righthand side is equal to $0$). Here we arranged the righthand side in the following order: 
\begin{itemize} 
\item $\mathcal{L}$ invading $\mathcal{S}$:  ``yes'' invading in \eqref{eq:small1}, and ``no'' invading in  \eqref{eq:small2}. 
\item $\mathcal{S}$ invading $\mathcal{L}$: ``yes'' invading in \eqref{eq:large1}, and ``no'' invading  in \eqref{eq:large2}. 
\end{itemize} 
Next, we perform some algebraic manipulations.  Before doing so, we stress that the three functions we are about to define, $S,L$ and $D$ are defined for all $(k,l)\in\Z\times \Z$. Consider the coefficients of $\nu(k,l)$ in 
on \eqref{eq:small1} and \eqref{eq:small2}. Their sum is 
$$ \frac{l}{m+n}\frac{k}{m} + \frac{n-l}{m+n}\frac{m-k}{m} = \frac{1}{m+n} \left ( n-l \frac{m-k}{m} - (n-l) \frac{k}{m}\right).$$
This leads to the following telescopic representation of \eqref{eq:small1} and \eqref{eq:small2}: 
\begin{align}
\label{eq:Spart}
\eqref{eq:small1} +\eqref{eq:small2} = \frac{n}{m+n} \nu(k,l)+ S(k,l) \end{align} 
 where 
 \begin{align}
 \label{eq:SKL} S(k,l) &= \frac{1}{(m+n)m} 
 \left ( \nu(k-1,l) l(m-k+1) - \nu(k,l)(n-l) k  \right)  \\
  & \quad+{ \frac{1}{(m+n)m} \left( \nu(k+1,l)(n-l) (k+1) -  \nu(k,l)l (m-k) \right) }.
  \end{align}
  Examining this expression, and defining 
  $$ D(k,l) = \frac{1}{(m+n)m}\left( \nu(k-1,l) l (m-k+1) -\nu(k,l) (n-l)k\right),$$
  we have that 
  \begin{equation}
  \label{eq:SDiff}
      S(k,l) = D(k,l)-D(k+1,l).
  \end{equation}
As an important consequence of this definition, we have that for all $l$, 
\begin{equation} 
\label{eq:Szero} \sum_{k=0}^m  S(k,l) = D(0,l)-D(m+1,0) = 0-0. 
\end{equation} 
Repeating the same analysis but with the roles of $k$ and $l$ changed, we arrive at 
\begin{align} 
\label{eq:Lpart}
\eqref{eq:large1} +\eqref{eq:large2}& = 
 \frac{m}{m+n} \nu(k,l)+ L(k,l) 
 \end{align} 
 where 
 \begin{align}
 \label{eq:Lone} L(k,l) &= \frac{k}{(m+n)n} 
 \left ( \nu(k,l-1) (n-l+1) - \nu(k,l) (n-l) \right) \\
 \label{eq:Ltwo}
  & + \frac{m-k}{(m+n)n} \left( \nu(k,l+1) (l+1) - \nu(k,l)l  \right), 
\end{align} 
Moreover, 
\begin{equation} 
\label{eq:Lzero}
\sum_{l=0}^n L(k,l)=0. 
\end{equation} 
Noting that the coefficients of $\nu(k,l)$ in \eqref{eq:Spart} \eqref{eq:Lpart} add up to $1$,   we can express the equation for the QSD in terms of $\mathcal{S}$ and $\mathcal{L}$ as follows: 
\begin{align}
\label{eq:clean_QSD1}
(\lambda -1 )\nu(k,l) &= S(k,l) + L(k,l)-{\bf 1}_\Delta(k,l)(S(0,0)+ L(0,0))
\end{align}
for all $(k,l) \in \{0,\dots,m\}\times\{0,\dots,n\}$. Note that the last term on the righthand side is necessary to guarantee that equation obtained also holds on the absorbing set. We have also used the symmetry of the model which gives $S(m,n)=S(0,0)$ and $L(m,n)=L(0,0)$. Indeed, 
$$(m+n)m S(m,n) \overset{\eqref{eq:SKL}}{=} \nu(m-1,n)n =\nu(1,0)n \overset{\eqref{eq:SKL}}{=}  (m+n)m S(0,0),$$
and an analogous derivation for $\mathcal{L}$. As a result, summing over all $k,l$ and using the fact that $\nu$ is a probability measure, we have 
$$(\lambda-1) = -2( S(0,0)+L(0,0)),$$
which in turn give us the equation 
\begin{equation}
\label{eq:clean_QSD}
(\lambda -1 )\nu(k,l) = S(k,l) + L(k,l)+\frac{\lambda-1}{2}{\bf 1}_\Delta(k,l).
\end{equation} 
\subsubsection{Switching to proportions} 
The next step is to switch from counting number of ``yes'' in $\mathcal{L}$ to the respective proportions. For each $n$   a number of $l$ ``yes'' corresponds to the proportion $l/n$. We will use the variable $x$ to represent these proportions. In \eqref{eq:QSD}  we introduced $\bar \nu$, the result of this change to $\nu$, and we similarly define $\bar S$, $\bar D$ and $\bar L$. That is 
$$\bar S(k,x) = S(k,nx){\bf 1}_{\{0,1/n,\dots,1\}}(x).$$ 
We will consider $\bar S$ and $\bar L$ as discrete measures on $\{0,\dots,m\}\times [0,1]$ and can therefore express  \eqref{eq:clean_QSD} as
\begin{equation} 
\label{eq:QSD_measure_rep} (\lambda-1) \int f d \bar \nu = \int f d \bar S + \int f d \bar L+ \frac{\lambda-1}{2}(f(0,0) + f(m,1)),
\end{equation} 
for every test function  $f:\{0,\dots,m\}\times [0,1] \to \R$. 
\subsubsection{Subsequential limits for $\mathcal{S}$}
\label{sec:Ssubseq}
From the construction of $D$, \eqref{eq:SDiff}, 
$$\int f d \bar S =\int f(k,x) d \bar D(k+1,x) - \int  f(k,x) d\bar D(k,x).$$ 
Take a function $f=f(k,x)$ which only depends on $x$. We will simply write it as $f=f(x)$. We have 
$$(\lambda-1) \int f (x) \bar \nu (k,dx) -\int f(x)    \bar L(k,dx) = \int f(x) \bar D(k+1,dx) - \int f(x)  \bar D(k,dx).$$ 
Take any subsequence along which $\bar \nu$ converges weakly to some limit $\bar \nu_\infty$. Denote the second marginal under $\bar \nu_\infty$ by $\bar\nu_{\infty,2}$. That is $\bar \nu_{\infty,2} (dx) = \sum _k\bar \nu_\infty(k,dx)$.  Without loss of generality we can also assume that the subprobability measures $\bar D$ also converge vaguely to some limiting subprobability measure $\bar D_\infty$.  Since $\lambda \to1$, it follows that 
$\bar D_\infty(k+1,dx)= \bar D_\infty(k,dx)$ for all $k$. By construction, this implies 
$(1-x) k \bar \nu_\infty(k,dx) =  x(m-k+1) \bar \nu_\infty (k,dx)$, and by solving the resulting recurrence relation we have
 $$\bar \nu_\infty(k,dx) = \binom{m}{k}x^k(1-x)^{-k} \bar \nu_\infty(0,dx).$$
 
Summing over $k$ and using the binomial formula, we have 
 $ \bar \nu_{\infty,2}(dx)  = (1-x)^{-m}\bar \nu_\infty(0,dx)$, and therefore we can rewrite the equation above as 
 $$ \bar \nu_\infty(k,dx) = \binom{m}{k} x^k (1-x)^{m-k} \bar \nu_{\infty,2}(dx),$$ 
 or, $\bar \nu_{\infty} (k|dx) \sim \mbox{Bin}(m,x)$. 
\subsubsection{Subsequential Limits for $\mathcal{L}$}
{\bf Rearrangement}\\
Let now $f=f(k,x)$ be only a function of the second variable. We will simply denote it by $f=f(x)$. We will assume that $f$ is also analytic on $\R$. Freezing $k$ and integrating $\bar L$ with respect to the second variable we have from \eqref{eq:Lone} and \eqref{eq:Ltwo} that 
\begin{align} (m+n) \int f(x) \bar  L (k,dx) & = \sum_{l'=0}^n \left ( f((l'+1)/n)  \frac{(n-l')}{n} k - f(l'/n) (m-k) \frac{l'}{n}  \right)\nu(k,l')\\
 & + \sum_{l'=0}^n \left ( f ((l'-1)/n) \frac{l'}{n} (m-k) -  f(l'/n)\frac{n-l'}{n} k \right)\nu(k,l')\\
 \label{eq:pretaylor1}
 & \overset{x=l'/n}{=} \int \left (  k (1-x)  f(x + \frac{1}{n})  - (m-k) x f(x) \right) \bar \nu(k,dx) \\
\label{eq:pretaylor2}
 & + \int \left ( (m-k) x f(x - \frac{1}{n}) - k (1-x) f(x) \right)\bar \nu(k,dx).  
 \end{align}
 {\bf Taylor Expansion}\\
 Using the Taylor expansion for $f$ at $(x\pm \frac{1}{n})$ up to the second order derivative, the error term is $o(n^{-3})$ uniformly in $k$ and $x$. After plugging  the Taylor expansion into \eqref{eq:pretaylor1} and \eqref{eq:pretaylor2}, we observe that the term $f(k,x)$ in the expansion cancels out and that we are left with  
\begin{align}
\nonumber \int f(x) d\bar L &= \frac{1}{(m+n)n} \int \left (k(1-x)-(m-k)x\right)  f'(x)  d \bar \nu \\
\nonumber
& + \frac{1}{2(m+n)n^2} \int \left ( k(1-x) + (m-k) x \right)f''(x) d \bar  \nu\\
\label{eq:taylor}
& + O(n^{-3}).
\end{align}
In addition, since $f$ is only a function of $x$, it follows from \eqref{eq:Szero} that  $\int f(x) d\bar S= 0$. Thus, plugging \eqref{eq:taylor} into \eqref{eq:QSD_measure_rep}, we have 
\begin{align} 
\nonumber
(\lambda-1) \int f (x) d \bar \nu &=  \frac{1}{(m+n)n} \int \left (k(1-x)-(m-k)x\right)  f'(x)  d \bar \nu  \\ 
\nonumber 
& +  \frac{1}{2(m+n)n^2} \int \left ( k(1-x) + (m-k) x \right)f''(x) d \bar  \nu\\
\label{eq:taylor_consq}
& + O(n^{-3})+ \frac{\lambda-1}{2} (f(0)+f(1)).
\end{align}
{\bf Equation for second marginal} \\
Since the left-hand side of \eqref{eq:taylor_consq} and the last three summands on the right-hand side are $O(n^{-3})$, the first integral on the right-hand side converges to zero. In other words, letting 
\begin{equation} c_n(f)  =  \int \left (k(1-x)-(m-k)x\right)  f'(x)  d \bar \nu,
\end{equation} 
then  $\lim_{n\to\infty} c_n(f)=0$. 
By replacing the function $f$ in \eqref{eq:taylor_consq} by  $f-c_nx$, we are left with 
\begin{align} 
(\lambda-1) \int f (x) -c_n x d \bar \nu 
& =  \frac{1}{2(m+n)n^2} \int \left ( k(1-x) + (m-k) x \right)f''(x) d \bar  \nu\\
& +(\lambda-1) (f(0)+f(1)-c_n)+ O(n^{-3}).
\end{align}
Multiplying both sides by $(n+m)n^2$, and taking a subsequential limit for $\bar \nu$ as in the previous step and using the fact that $\lambda -1 \sim -2m n^{-3}$ we obtain 
\begin{equation}
    \label{eq:ototo}
-2m \int f \bar d\nu_{\infty,2}=\frac 12 \int \left ( k(1-x) + (m-k) x \right)f''(x) d \bar  \nu_{\infty}- 2m (f(0)+f(1))
\end{equation}
Moreover, using the result of Section \ref{sec:Ssubseq}, we have 
$ \int k \bar \nu_\infty(dk|x) = mx$, the expectation of $\mbox{Bin}(m,x)$, and so the first expression on  right-hand side of \eqref{eq:ototo} becomes 
$$ 2m \int (x-x^2) f''(x) d \bar \nu_{\infty,2}(x).$$ 
Plugging this into \eqref{eq:ototo} and dividing by $m$ yields 
\begin{equation} 
\label{eq:steins} 
\int x(1-x) f''(x)  +2f(x) d \bar \nu_{\infty,2}(x)=f(0)+f(1).
\end{equation} 
\subsubsection{Identifying the measure} 
Now fix any  $F$ which is supported on a closed sub-interval of $(0,1)$. The differential equation 
\begin{equation} 
\label{eq:ODE} x(1-x)f'' (x) + 2f (x) = F(x)
\end{equation} 
has a  continuous solution on $[0,1]$.  It is noted that all solutions are bounded and continuous because two linearly independent solutions to the homogeneous equation are such. In order to avoid disruption to the flow we will prove this standard result in ordinary differential equations in the Appendix.  As a result, we have that $F$ is integrable with respect to $\bar \nu_{\infty}$ and that 
$$\int F (x) d \bar\nu_{\infty,2} = f(0)+f(1).$$ 
Also, 
$$ \int_0^1 F (x) dx  =\int_0^1 x(1-x) f''(x) + 2f(x) dx$$ 
Integrating the right-hand side by parts  we obtain 
\begin{equation} 
\begin{array}{lll} 
\int_0^1 x f''(x)dx &- \int_0^1 x^2 f''(x) dx &+ 2 \int_0^1 f(x) dx\\
 = f'(1)-(f(1)-f(0)) &- (f'(1) - 2 \int_0^1 x f'(x) dx)  &+ 2 \int_0^1 f(x) dx \\ 
 = f(0)-f(1) &+ 2f(1)-\int_0^1 f(x) dx  &+ 2 \int_0^1 f(x) dx \\ 
= f(0)+f(1)
\end{array} 
\end{equation} 
Therefore, we have proved that 
$$\int F(x) d \bar\nu_{\infty,2}  =f(0)+f(1) = \int_0^1 F(x) dx.$$ 
It therefore follows that $\bar\nu_{\infty,2}$ is uniform on $[0,1]$. 
\section{Appendix}
\subsection{Solution to ODE \eqref{eq:ODE}} 
We show that \eqref{eq:ODE} has a continuous and bounded solution on $[0,1]$. To this end we begin with a change of variables that will transform the equation to something more familiar. Let $t= 2x-1$. Since $x \in [0,1]$, $t \in [-1,1]$. Now let $v(t) = f(x) = f(\frac{t+1}{2})$ and let $G(t) = F(x) = G(\frac{t+1}{2})$. Now $\frac{dv}{dt} = \frac{d f}{dx}\frac{dx} {dt} =\frac{1}{2} f'(x)$, and so $v''(t) = \frac 14  f''(x)$. This gives 

$$x(1-x) f''(x) + 2 f(x) = \frac{1-t^2}{4} 4 v''(t) + 2 v(t),$$ 
Or 
\begin{equation} 
\label{eq:ODEt} (1-t^2) v''(t) + 2 v(t)=G(t)
\end{equation} 

Two linearly independent solutions on $(-1,1)$ are 
$$ v_1 (t) = t^2-1,\quad v_2 (t) = 2t + (t^2-1)  \ln \frac{1-t}{1+t}.$$ 
Both solutions are continuous on $[0,1]$. 

A particular solution to can be found using variation of parameters, that is the solution is of the form $v(t) = c_1 (t) v_1(t) + c_2(t) v_2(t)$. Any other solution differs from it by a linear combination of $v_1$ and $v_2$. 

Let $W$ be the Wronskian for $v_1,v_2$. That is 

$$W(t)=\begin{pmatrix} v_1 & v_2 \\ v_1' & v_2'\end{pmatrix}(t).$$ 

Then the coefficient functions of $c_1,c_2$ of the particular solution are obtained through the equation 

$$\begin{pmatrix} c_1 \\ c_2 \end{pmatrix}(t) = \int_0^t G(s) W^{-1}(s) \begin{pmatrix}0 \\ 1 \end{pmatrix}ds. $$ 
Since $W^{-1}$ is continuous on $(-1,1)$ and $G$ is compactly supported, $c_1$ and $c_2$ are bounded, continuous functions on $[-1,1]$. This completes the proof. 

\subsection{Numerical Approximation} 
\subsubsection{Survival Rates} 
\begin{figure}[hb]
\centering
\resizebox{\textwidth}{!}{
\includegraphics{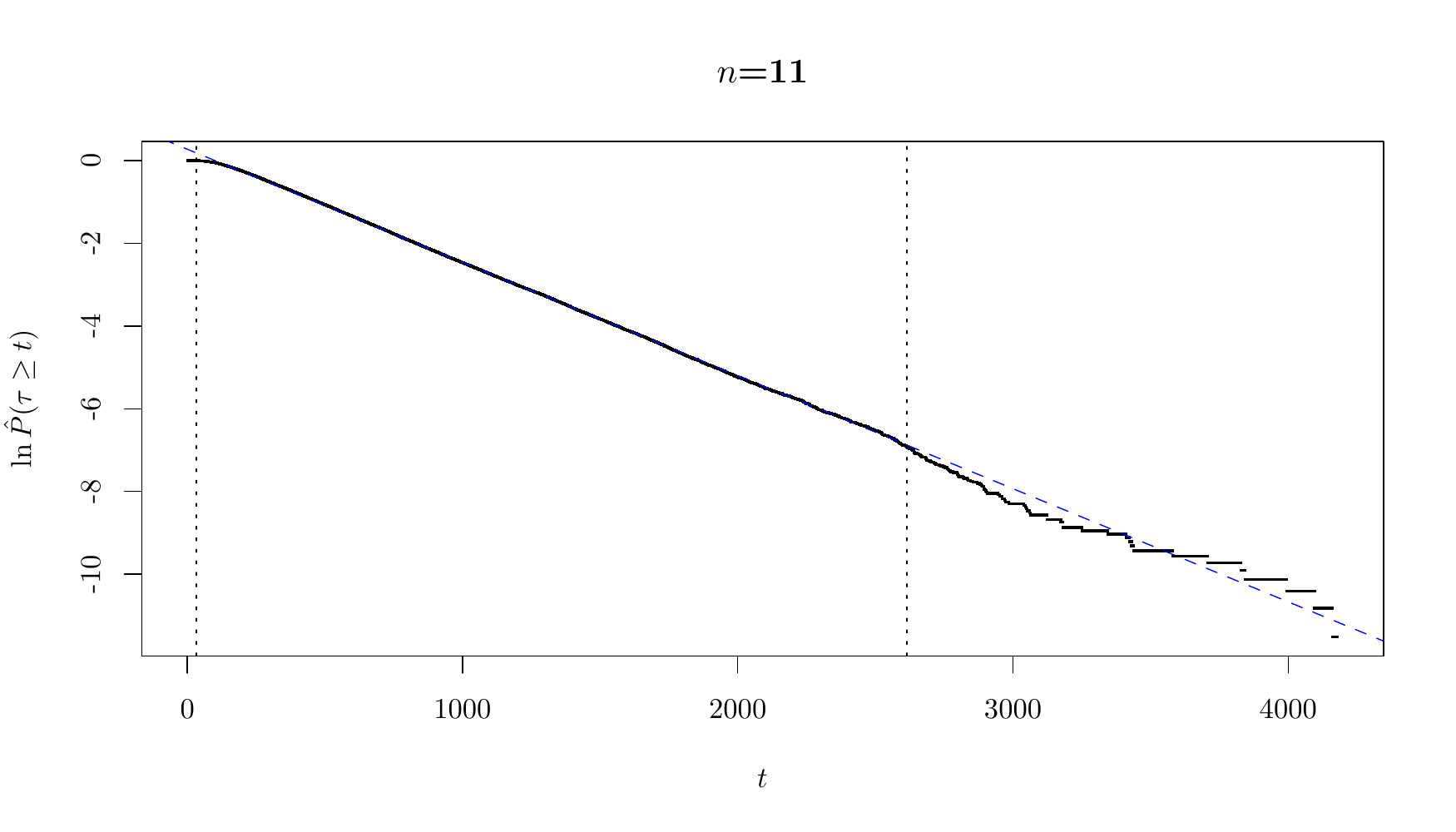}
\includegraphics{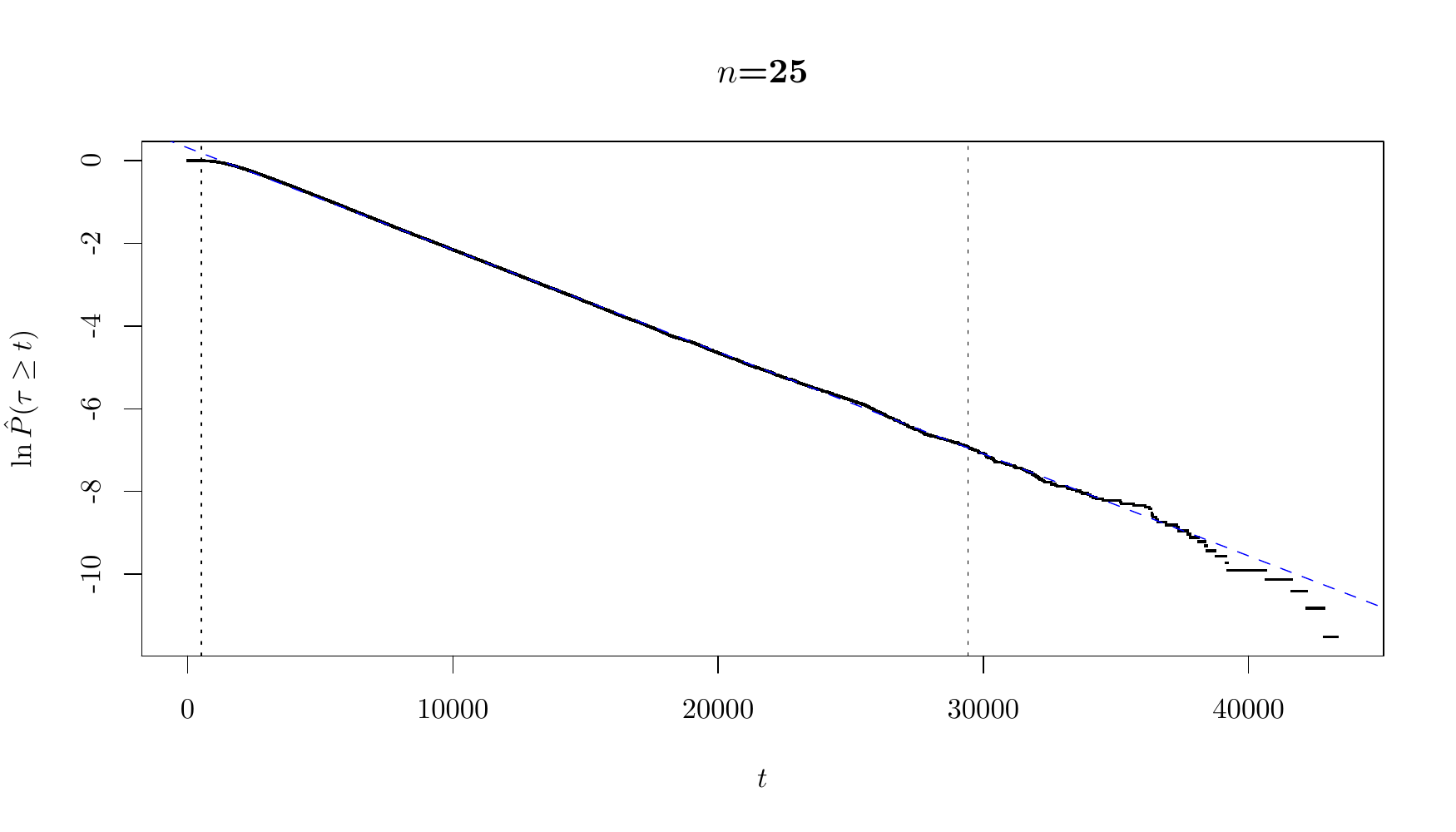}
}
\caption{Log-plot of the tail of the empirical distribution $t\to \hat P(\tau>t)$ for the consensus time for the Invasion Model on $K_{2,n}$ for $n=11,25$, in black. Each plot was obtained from $10^5$ simulated paths, with an initial configuration of $\lceil n/2\rceil$  and $1$ ``yes'' opinions in  ${\cal L}$ and ${\cal S}$, respectively,  and $\hat P(\tau\ge t)$ is the proportion of paths whose consensus time is larger or equal to $t$. We performed a linear regression on these data, omitting the lowest and highest $0.1\%$ percentiles for each empirical distribution, marked by dashed vertical lines.  The resulting regression line is in blue and its slope gives an estimate of $\ln (1-\lambda_{OD}(\rho_I))$. The results of this procedure applied to $n=11,13,15,\dots,25$ are presented in Figure \ref{fig:loglogplot}.}
\label{fig:logplots}
\end{figure}

\begin{figure}[h]
\centering
\resizebox{\textwidth}{!}{
\includegraphics{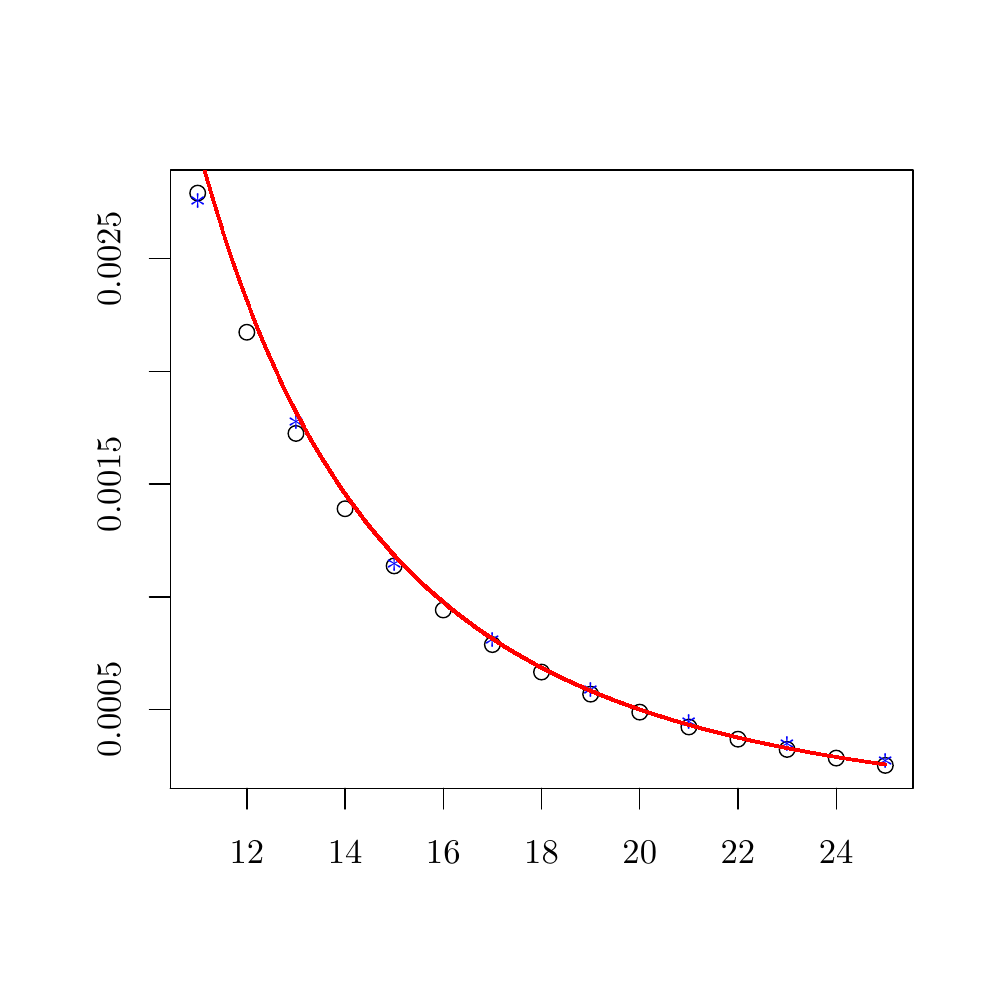}
\includegraphics{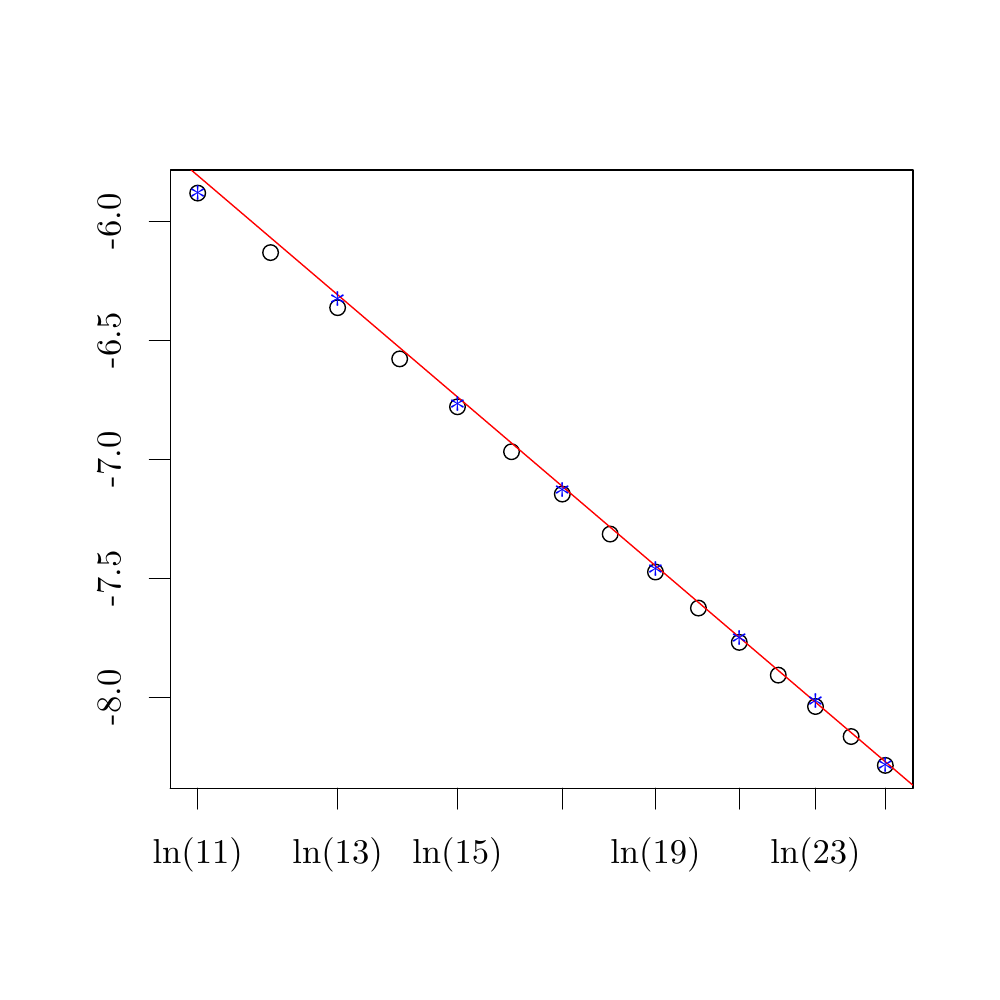}}
\caption{Graphs for $1-\lambda_{OD}(\rho_I)$ for $K_{2,n}$, as a function of $n$, left, and the corresponding $\log-\log$ plot, right. The  circles represent values obtained by a numerical calculation of the eigenvalue. The asterisks (only for odd values of $n$) are estimated values obtained from the linear regression discussed  in Figure \ref{fig:logplots}. The line represents the asymptotic curve $\frac{4}{n^3}$ as $n\to\infty$ from equation \eqref{eq:lambdaI}.}
\label{fig:loglogplot}
\end{figure}

The  sharp estimate for the survival rate  $\lambda_{OD}(\rho_I)$ in  Proposition \ref{prop:lambdacmc} is a key component  in the proof of our main result Theorem \ref{thm:main}. Our work began by trying to identify the dependence of the survival rate $\lambda_{OD}(\rho_I)$ for the Invasion Model on $K_{m,n}$ on $n$ and $m$ empirically.   This was done through standard Monte-Carlo simulation for the tail of the absorption time $\tau$.  As the tail of $\tau$ decays geometrically with rate $\lambda_{OD}$, \eqref{eq:log_tail}, one expects the log-plot of the empirical tail function for $\tau$ to be asymptotically linear with slope $\ln \lambda_{OD}(\rho_I)$. Thus, an estimate for  $\lambda_{OD}(\rho_I)$ can be obtained by linear regression. Figure \ref{fig:logplots} shows the results of some of these simulations, including the resulting regression curves, whose slopes are very close to the asymptotic formula from \eqref{eq:lambdaI} from Proposition \ref{prop:lambdacmc}. Using an array of such  simulations with different values of  $m$ and $n$ we were able to identify an empirical dependence of the form  $1-\lambda_{OD}(\rho_I)\sim \frac{2m}{n^3}$. Figure \ref{fig:loglogplot}  does that. In the figure we used  both the estimated regression values, numerically computed values,  and values obtained from the asymptotic formula \eqref{eq:lambdaI} in Proposition \ref{prop:lambdacmc}. 

\subsubsection{The QSD}
Numerical approximation of QSDs is often performed through simulation or through numerical solution of the corresponding eigenvector equation. As QSDs naturally  appear as limits (QLDs) of distributions conditioned on events with eventually vanishing probabilities, sampling  QSDs through simulations is non-trivial. The naive approach obtained by sampling paths and keeping only those not absorbed by some long time is extremely inefficient as nearly all sampled data is eventually discarded. Numerical solutions to the eigenvector equation is efficient when dealing with small systems but become very cumbersome when the systems are very large, also because it is often also the case where the effects of absorption are finer and demand increasing precision.  
\subsubsection{Survival Rates}
In our work we implemented a combination of the latter method for small values of $m$ and $n$, and a simulation scheme introduced in \cite{aldous} (further developed and generalized in \cite{benaim}) which is easy to implement, and which we found to be very efficient and precise. The idea is the following: run one simulated copy of the Invasion Model until consensus, but not moving to the consensus state. Instead, start the process afresh from a distribution  on the non-absorbing states. The distribution used is the empirical distribution up to that time: the probability of starting from any state is the  proportion of time the process spent at that state. This procedure is repeated indefinitely. It turns out that  the empirical distribution converges a.s. to the QSD. Through this method we were able to obtain a good picture of the structure of the QSD for large values of $n$, which eventually helped us in stating our main result. These simulations were also very helpful in earlier steps when studying the dependence of the survival rates on $m$ and $n$. 
\subsubsection{Other directions}
Though here our primary focus was the QSD for the Invasion Model,  the work led to other directions, even within the context of dynamics on the complete bipartite graph. We list three here. 

\begin{figure}
    \resizebox{\textwidth}{!}{
  \includegraphics{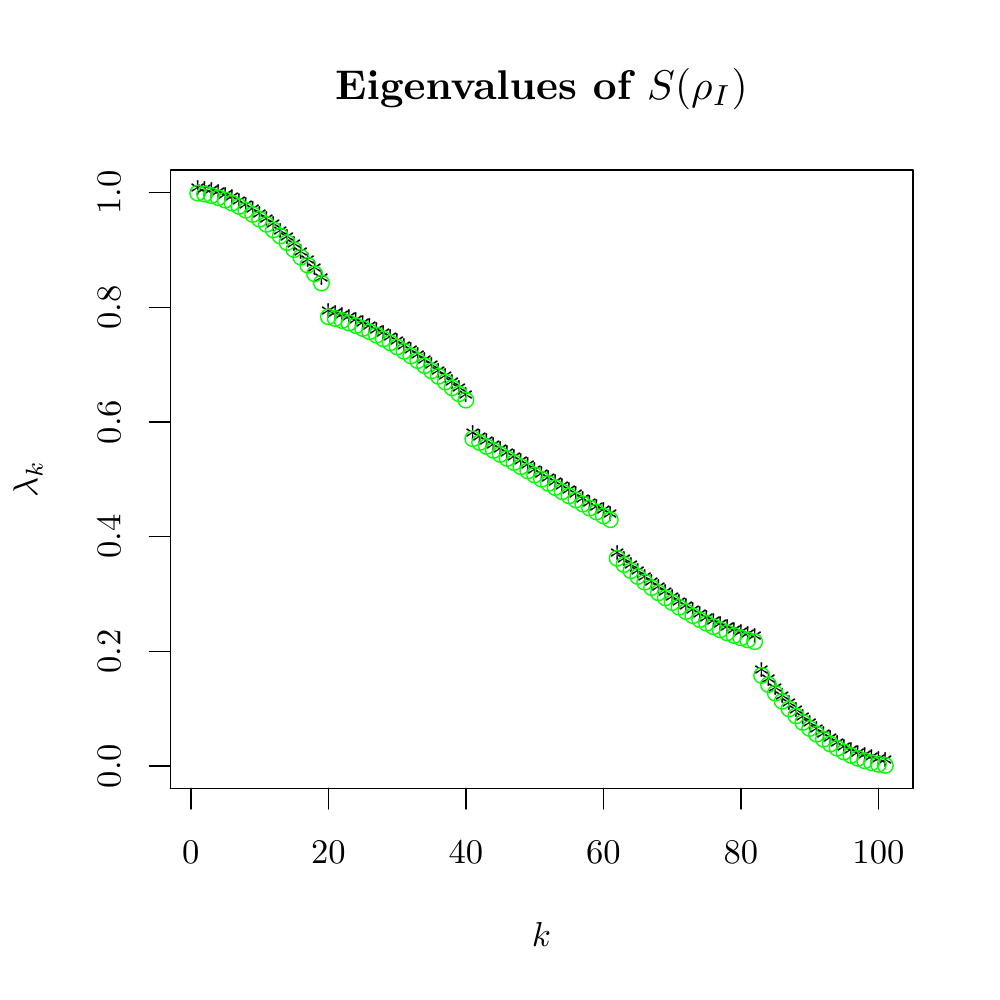}}
   \caption{Eigenvalues for $S(\rho_I)$ on $K_{4,20}$, represented by black $*$'s and their reflections with respect to the  middle point between the minimum and the maximum. The rows and columns in $S(\rho_I)$ corresponding to the inaccessible states $(4,0)$ and $(0,20)$ were removed from $S(\rho_I)$.}
    \label{fig:spectrum}
\end{figure}

\begin{enumerate}
\item One concerns models obtained by interpolating the Voter and Invasion model. That is, at each step, we sample a directed edge according to the Voter Model with some probability $p$ and according to the Invasion Model with probability $1-p$. 
\item What happens in the ``intermediate'' stage? In the Invasion Model the small group converges rather quickly to the limiting distribution of conditioned on the density of ``yes" large group, because the members of which are very frequently invaded by the large group, while the large group varies very slowly. Similarly, in the Voter Model, the small group reaches consensus rather quickly. What are the time scales involved? Is meta-stability exhibited, namely, for a long-stretch of time the system  (conditioned to not to reach consensus) is near a  stationary/stable state, different from the QSD. Do the dynamics exhibit a cut-off phenomenon where the distribution of the process abruptly switches from being orthogonal to the QSD to being close to the QSD? 
\item The spectrum of $S(\rho_I)$ for the Invasion Model calculated numerically appears to exhibits some structure.   Figure \ref{fig:spectrum} shows the eigenvalues  of $S(\rho_I)$  on $K_{m=4,n=20}$, after the two inaccessible states (with the inaccessible states $(m,0)$ and $(0,n)$ were removed. As we also eliminated the absorbing states $(0,0)$ and $(m,n)$, there are $(m+1)*(n+1)-4= 101$ remaining states. The following observations appeared for other values of $m$ and $n$ we examined. 
\begin{itemize} 
\item The number of distinct eigenvalues is equal to the number of states, and all are real. Denote them by  $\lambda_{OD}(\rho_I)=\lambda_1>\lambda_2>\dots >\lambda_{(m+1)*(n+1)-4}$, each appearing in the figure as an $*$. 
\item The eigenvalues are split into $m+1$ groups separated by gaps. 
\item The eigenvalues exhibit a symmetry. The numerical calculations appear to suggest that the set is symmetric with respect to the middle point $(\lambda_1 + \lambda_{(m+1)*(n+1)-4})/2$. The eigenvalues, reflected about the middle, and arranged in decreasing order appear in green circles. 
\end{itemize}
\end{enumerate} 
\clearpage
\bibliographystyle{amsalpha} 
\bibliography{main.bib}
\end{document}